\UseRawInputEncoding
\documentclass[12pt]{article}
\usepackage{amsmath,amssymb,amsthm, amsfonts}
\usepackage{hyperref}
\usepackage{graphicx}
\usepackage{url}
\usepackage{dsfont} 
\usepackage{xcolor}
\definecolor{red}{rgb}{1,0,0}

\definecolor{blue}{rgb}{0,0,1}

\definecolor{green}{rgb}{0,.6,0}

\usepackage{float}
\usepackage{tikz}

\newtheorem{theorem}{Theorem}[section]

\newtheorem{lemma}[theorem]{Lemma}

\newtheorem{observation}[theorem]{Observation}

\newtheorem{problem}[theorem]{Problem}

\newtheorem{construction}[theorem]{Construction}
\newtheorem{remark}[theorem]{Remark}

\newcommand{\sat}{\mathrm{sat}}
\newcommand{\BF}{\text{Berge-$F$}}
\newcommand{\mC}{\mathcal{C}}
\newcommand{\mS}{\mathcal{S}}
\newcommand{\mF}{\mathcal{F}}
\newcommand{\mH}{\mathcal{H}}
\newcommand{\mG}{\mathcal{G}}
\newcommand{\mK}{\mathcal{K}}

\title{Saturation Numbers for Berge Cliques}
\author{Sean English\thanks{Dept.~of Mathematics and Statistics, University of North Carolina Wilmington (\texttt{EnglishS@uncw.edu})}\and 
	J\"urgen Kritschgau\thanks{Dept.~of Mathematics, Carnegie Mellon University (\texttt{jkritsch@andrew.cmu.edu}). Research is supported by NSF grant DMS-1839918.}\and Mina Nahvi\thanks{Dept.~of Mathematics, University of Illinois Urbana-Champaign (\texttt{mnahvi2@illinois.edu})} \and Elizabeth Sprangel\thanks{Dept.~of Mathematics, Iowa State University (\texttt{sprangel@iastate.edu})}}

\date{}

\begin{document}

	\maketitle
	
	\begin{abstract}
		Let $F$ be a graph and $\mH$ be a hypergraph, both embedded on the same vertex set. We say $\mH$ is a Berge-$F$ if there exists a bijection $\phi:E(F)\to E(\mH)$ such that $e\subseteq \phi(e)$ for all $e\in E(F)$. We say $\mH$ is Berge-$F$-saturated if $\mH$ does not contain any Berge-$F$, but adding any missing edge to $\mH$ creates a copy of a  Berge-$F$. 
		The saturation number $\sat_k(n,\text{Berge-}F)$ is the least number of edges in a Berge-$F$-saturated $k$-uniform hypergraph on $n$ vertices. 
		We show
		\[
		\sat_k(n,\text{Berge-}K_\ell)\sim \frac{\ell-2}{k-1}n,
		\]
		for all $k,\ell\geq 3$. Furthermore, we provide some sufficient conditions to imply that $\sat_k(n,\text{Berge-}F)=O(n)$ for general graphs $F$.
	\end{abstract}
	
	\section{Introduction}
	
	Extremal graph theory is concerned with maximizing or minimizing some parameter over a restricted class of graphs. Let $\mG$ and $\mF$ be $k$-uniform hypergraphs. We say that $\mG$ is $\mF$-saturated if $\mG$ does not contain a copy of $\mF$ but $\mG+e$ does for any $e\in E(\overline{\mG})$.
	The most well-studied problem in extremal graph theory is the Tur\'an problem, which asks for the maximum number of edges in a $\mF$-free hypergraph $\mG$ on $n$ vertices. This maximum is known as the \textbf{extremal number} or \textbf{Tur\'an number} of $\mathcal{F}$, and is denoted $\mathrm{ex}_k(n,\mathcal{F})$.
	Any $\mF$-free hypergraph $\mG$ with $\mathrm{ex}_k(|V(\mG)|,\mathcal{F})$ edges must necessarily be $\mathcal{F}$-saturated, so we can write
	\[
	\mathrm{ex}_k(n,\mathcal{F})=\max\{|E(\mG)|: |V(\mG)|=n, \mG\text{ is }\mF\text{-saturated}\}.
	\]
	On the flipside, the \textbf{saturation number} of $\mF$, denoted $\sat_k(n,\mathcal{F})$, is the least number of edges in a $\mF$-saturated graph on $n$ vertices, or
	\[
	\sat_k(n,\mathcal{F})=\min\{|E(\mG)|:|V(\mG)|=n,\mG\text{ is }\mF\text{-saturated}\}.
	\]
	Saturation was first introduced by Erd\H{o}s, Hajnal and Moon~\cite{EHM64} for graphs, and then generalized for hypergraphs by Bollob\'as~\cite{B65} who showed that
	\begin{equation}\label{equation saturation for complete hypergraph}
		\sat_k(n,\mK^{(k)}_\ell)=\binom{n}{k}-\binom{n-\ell+k}{k},
	\end{equation}
	where $\mK^{(k)}_\ell$ denotes the complete $k$-uniform hypergraph on $\ell$ vertices. Since these seminal results, much work has been done on the saturation function and many generalizations have been studied. For a dynamic survey on saturation numbers, see~\cite{FFS11}.
	
	In this work, we are interested in the saturation function for Berge hypergraphs, which are a generalization of Berge paths and Berge cycles introduced by Gerbner and Palmer~\cite{GP2017}.
	Given a graph $F$ and a hypergraph $\mH$ embedded on the same vertex set, we say that $\mH$ is a \textbf{Berge-$F$} if there is a way to embed $F$ and $\mH$ on the same vertex set such that there exists a bijection $\phi:E(F)\to E(\mH)$ that has $e\subseteq \phi(e)$ for all $e\in E(F)$. 
	We note that many non-isomorphic hypergraphs may be a Berge-$F$ and a hypergraph may be a Berge copy of many non-isomorphic graphs. We will write $\sat_k(n,\text{Berge-}F)$ for the least number of edges in a Berge-$F$-saturated $k$-uniform hypergraph on $n$ vertices.
	
	Saturation numbers for Berge hypergraphs were first studied by the first author and others in \cite{EGKMS2017}, where some results on the saturation function for Berge paths, matchings, cycles, and cliques are given. Since the seminal work on saturation for Berge hypergraphs, the topic has gotten significant attention (see \cite{AE2019, AW2018, GPTV2022, LHK2020} for some of the results on the topic). Prior to this work, the following result was the only known result on saturation for Berge cliques.
	
	\begin{theorem}[\cite{EGKMS2017}]\label{theorem berge triangles}
		For all $k\geq 3$ and $n\geq k+1$,
		\[
		\sat_k(n,\text{Berge-}K_3)=\left\lceil\frac{n-1}{k-1}\right\rceil.
		\]
	\end{theorem}
	
	Our main theorem determines the asymptotics of Berge-$K_\ell$ for all fixed clique sizes $\ell$ and uniformities $k$.
	
	\begin{theorem}\label{theorem main}
		For $\ell\geq 3$ and $k\geq 3$, 
		\[
		\sat_k(n,\text{Berge-}K_\ell)\sim\frac{\ell-2}{k-1}n.
		\]
	\end{theorem}
	
	The case $\ell=3$ is covered by Theorem~\ref{theorem berge triangles}. When $\ell\geq 4$, the lower bound in Theorem~\ref{theorem main} follows from Theorem~\ref{theorem lower bound}, while the upper bound follows from Theorem~\ref{theorem upper bound}.
	
	In addition to the main result, we also study the linearity of Berge saturation for general graphs. For ($2$-uniform) graphs $F$, K\'aszonyi and Tuza~\cite{KT86} showed that $\sat_2(n,F)=O(n)$, while Pikhurko~\cite{P99} showed that $\sat_k(n,\mF)=O(n^{k-1})$ for $k$-uniform hypergraphs $\mF$, and this result is best-possible, as seen for example in equation~\eqref{equation saturation for complete hypergraph} stating the result from \cite{B65}. In~\cite{EGKMS2017} it was conjectured that $\sat_k(n,\text{Berge-}F)=O(n)$, suggesting that the saturation function for Berge hypergraphs should grow more like graph saturation than $k$-uniform hypergraph saturation. This conjecture was confirmed for uniformities $3\leq k\leq 5$ in~\cite{EGMT2018}, but is still open in general.
	
	We prove two results which show that many  graphs have Berge saturation numbers which grow at most linearly. The first theorem deals with graphs with large minimum degree.
	
	\begin{theorem}\label{theorem linearity delta bound}
		If $\delta(F)>\frac{|V(F)|-\alpha(F)}{2}$ and $k\geq 3$, then 
		\[
		\sat_k(n,\text{Berge-}F)=O(n).
		\]
	\end{theorem}
	
	Our final result concerns graphs with large girth. Recall that the \textbf{girth} of a graph is the length of the shortest cycle (which we will denote by $g(G)$), and the \textbf{vertex feedback number} is the least number of vertices necessary to delete which leaves an acyclic graph (which we will denote by $f(G)$).
	
	\begin{theorem}\label{theorem linearity girth bound}
		Let $F$ be a 
		graph with girth $g$ and vertex feedback number $f$. If $g> f$ and $k\geq 3$, then 
		\[
		\sat_k(n,\BF)=O(n).
		\] 
	\end{theorem}
	
	\subsection{Definitions and organization}
	Let $F$ be a graph and $\mH$ be a $k$-uniform Berge-$F$ embedded on the same vertex set, and let $\phi:E(F)\to E(\mH)$ be a bijection such that $e\subseteq \phi(e)$ for all $e\in E(F)$. We call $\phi$ the \textbf{Berge edge map}.
	When $F$ and $\mH$ are embedded in such a way that there exists a Berge edge map, we say that $F$ is a \textbf{Berge-$F$ witness}. 
	When $F$ is a Berge-$F$ witness for $\mH$, the vertices in $V(F)$ are called \textbf{core vertices} of the Berge-$F$ hypergraph $\mH$. 
	Given a hypergraph $\mH$ and a set $e\subseteq 2^{V(\mH)}$ such that $e\notin E(\mH)$, we say $\mH+e$ contains a \textbf{new} Berge-$F$ if $\,H+e$ contains a Berge-$F$ that uses $e$. 
	
	Uniform hypergraphs are of primary concern in this paper, but in order to simplify proofs, we find it useful to occasionally deal with non-uniform hypergraphs, usually a hypergraph where all but one edge has $k$ vertices, and the one other edge has $2$ vertices. We note that the definition of a Berge-$F$ does not depend on $\mH$ being uniform. In particular, we will say a pair $uv\subseteq V(\mH)$ is \textbf{$\ell$-good} if $\mH+uv$ contains a new Berge-$K_\ell$. Since we occasionally speak of non-uniform hypergraphs, we note here that if we refer to the complement $\overline{\mH}$ of a $k$-uniform hypergraph, this is assumed to be the $k$-uniform complement.
	
	In a $k$-uniform hypergraph $\mH$, a \textbf{loose path} of length $2$ is a pair of edges that intersect in exactly one vertex. The single vertex of degree $2$ in the loose path of length $2$ will be called a \textbf{hinge vertex}. Given two vertices $u,v\in V(\mathcal{H})$, we will write $v\preceq u$ if $\{e\in E(\mH)\mid v\in e\}\subseteq \{e\in E(\mH)\mid u\in e\}$. We note that $\preceq$ is reflexive and transitive, but in general $\preceq$ may not be a partial order as it may not be anti-symmetric. 
	
	We note that if $F$ is a non-empty graph with isolated vertices and $F'$ is the subgraph of $F$ induced by $E(F)$, then for all $n$ large enough, $\sat_k(n,\text{Berge-}F)=\sat_k(n,\text{Berge-}F')$, so throughout the paper we will silently assume that no graphs have isolated vertices. All asymptotics are with respect to $n\to\infty$, with all other parameters assumed to be constant unless specifically stated otherwise. All logarithms written as $\log n$ are in base $2$. The \textbf{complete join} of two $2$-graphs, $F$ and $G$, denoted $F\vee G$, is the graph whose vertex set is the disjoint union $V(F\vee G)=V(F)\cup V(G)$, and the edge set 
 \[
 E(F\vee E)=E(F)\cup E(G)\cup \{xy\mid x\in V(F), y\in V(G)\}.
 \]
	
	The rest of the paper is organized as follows. In Sections~\ref{section lower bound} and \ref{section upper bound}, we prove the lower bound and upper bound for Theorem~\ref{theorem main}, respectively.
	In Section~\ref{section linearity}, we prove Theorems \ref{theorem linearity delta bound} and \ref{theorem linearity girth bound}. Finally, in Section~\ref{section conclusion}, we briefly discuss the open problem of determining the exact values for the saturation numbers for Berge-$K_4$.
	
	\section{Lower bound for Berge-\texorpdfstring{$K_\ell$}{2} saturation}\label{section lower bound}
	
	We present here an asymptotic lower bound for $\sat_k(n,\text{Berge-}K_\ell)$. 
	Fortunately, the same argument works for all $k\geq 2$ and $\ell\geq 3$.
	
	\begin{theorem}\label{theorem lower bound}
		For any $k\geq 2$ and $\ell\geq 3$,
		\[
		\sat_k(n,\text{Berge-}K_\ell)\geq (1+o(1))\frac{\ell-2}{k-1}n.
		\]
	\end{theorem}

	\begin{proof}
		Let $\mH$ be a $k$-uniform Berge-$K_\ell$-saturated hypergraph. Partition the vertex set $V(\mH)=X\cup A\cup B$, where
		\[
		X=\{v\in V(\mH)\mid d(v)\geq \log^2 n\},
		\]
		$A\subseteq V(\mH)\setminus X$ such that $v\in A$ if and only if $v$ is contained in at least $\ell-2$ edges that intersect $X$, and $B=V(\mH)\setminus (X\cup A)$.
		Note that if $|X|> n/\log n$, then by counting degrees, $|E(\mH)|>\frac{|X|\log^2 n}k=\omega(n)\geq   (1+o(1))\frac{\ell-2}{k-1}n$, so we are done unless $|X|\leq n/\log n=o(n)$. 
		
		We will show that $|B|=o(n)$ as well. Indeed, first note that every non-edge $f\subseteq B$ contains a pair $u,v\in f$ that is connected by $\ell-2$ Berge paths of length $2$ since adding $f$ to $\mH$ creates a Berge-$K_\ell$.
		Since $u\not\in A$, at least one of these Berge paths must contain a hinge vertex in $A\cup B$. Each vertex $a\in A\cup B$ can play the role of this hinge vertex for at most $\binom{d(a)}2(k-1)^2\binom{|B|}{k-2}$ size $k$ non-edges $f\subseteq B$ since we can choose two edges in the Berge path of length $2$ in $\binom{d(a)}2$ and then the vertices $u$ and $v$ in at most $(k-1)$ ways each, and finally the remaining $k-2$ vertices in $f$ in $\binom{|B|}{k-2}$ ways. So, if $p$ is the total number of $k$-sets $f\subseteq B$ that are not edges of $\mH$, then
		\begin{equation}\label{equation lower bound upper bound on p}
			p\leq \sum_{a\in A\cup B}\binom{d(a)}2(k-1)^2\binom{|B|}{k-2}\leq k^2\binom{|B|}{k-2}\sum_{a\in A\cup B}\binom{\log^2 n}2\leq k^2\binom{|B|}{k-2}\frac{n\log^4 n}2.
		\end{equation}
		On the other hand, since every vertex $b\in B$ has degree less than $\log^2 n$, by a degree count, $B$ contains at most $\frac{|B|\log^2 n}k$ edges, and thus
		\begin{equation}\label{equation lower bound lower bound on p}
			\binom{|B|}{k}-\frac{|B|\log^2 n}k\leq p.
		\end{equation}
		If $|B|\geq n/\log n$, then the left side of \eqref{equation lower bound lower bound on p} is greater than $\frac{1}{2}\binom{|B|}{k}$. Comparing this to the right side of \eqref{equation lower bound upper bound on p} with some rearranging, we get that
		\[
		\frac{(|B|-k+2)(|B|-k+1)}{k^3(k-1)}\leq n\log^4 n,
		\]
		which is a contradiction for $|B|\geq n/\log n$. Thus, $|B|=o(n)$ as claimed.
		
		Thus, since $|X|,|B|=o(n)$, we must have $|A|=(1+o(1))n$. Now, let us count the number of edges of $\mH$ that contain at least one vertex in $A$ and at least one vertex in $X$. Each such edge contains at most $k-1$ vertices in $A$, and each vertex in $A$ is in at least $\ell-2$ such edges by definition of $A$, so the total number of edges containing at least one vertex from $X$ and at least one from $A$ is at least
		\[
		\frac{\ell-2}{k-1}|A|\leq |E(\mH)|.
		\]
		Since $|A|=(1+o(1))n$, the theorem holds.
	\end{proof}
	
	\section{Upper bound}\label{section upper bound}
	
	In this section, we provide Berge-$K_\ell$-saturated constructions for all uniformities $k\geq 3$ and all clique sizes $\ell\geq 4$.
	The work is divided into three subsections.
	In Section~\ref{subsection upper bound tools}, we state and prove a few simple observations which will be useful in the later sections. 
	In Section~\ref{subsection upper bound small graphs}, we construct specific small $k$-uniform Berge-$K_\ell$-saturated hypergraphs which also have the properties that every pair is $\ell$-good, and that every set of $\ell-1$ vertices contains a Berge clique.
	Due to some technical constraints when $\ell=4$, we provide one construction for $\ell=4$ and one for $\ell\geq 5$.
	In Section~\ref{subsection upper bound final construction}, we use the small hypergraphs constructed in Section~\ref{subsection upper bound small graphs} to find a Berge-$K_\ell$ saturated hypergraph with few edges on $n$ vertices for all large $n$.
	
	\subsection{Upper bound tools}\label{subsection upper bound tools}
	
	We present here two simple observations which will help show that our constructions are Berge-$K_{\ell}$-saturated.
	
	\begin{observation}\label{observation vertex encompass adding an edge}
		Let $\mH$ be a hypergraph and let $u,v\in V(\mH)$ be such that $v\preceq u$. 
		Let $e\in 2^{V(\mH)}$ be such that $v\in e$.
		Let
		\[
		e'=\begin{cases}e&\text{ if }u\in e,\\
			(e\setminus \{v\})\cup \{u\}&\text{ if }u\not\in e.
		\end{cases}
		\]
		For any graph $F$, if $\mH+e$ contains a new Berge-$F$ in which $u$ is not core, then $\mH+e'$ contains a new Berge-$F$. Furthermore, there is a new Berge-$F$ in $\mH+e$ and a new Berge-$F$ in $\mH+e'$ such that the two Berge-$F$'s have the same core vertices except possibly with $u$ as a core vertex instead of $v$.
	\end{observation}
	
	\begin{proof}
		Assume $H+e$ contains a Berge-$F$, call it $\mathcal{F}$, that uses the edge $e$ and in which $u$ is not core.
		If $v$ is not core in $\mathcal{F}$, then $(\mathcal{F}-e)+e'$ is a Berge-$F$ with the same witness as $\mathcal{F}$.
		If $v$ is core, then since every edge containing $v$ also contains $u$, $(\mathcal{F}-e)+e'$ is a Berge-$F$ with the same core vertices as $\mathcal{F}$, except with $u$ in place of $v$. 
	\end{proof}
	
	\begin{observation}\label{observation encompass F to F}
		Let $\mH$ be a hypergraph and let $u,v\in V(H)$ be such that $v\preceq u$.
		If $\mH$ contains a Berge-$F$ whose core vertices include $v$ and not $u$, then $\mH$ also contains a Berge-$F$ with the same core vertices, except with $u$ in place of $v$.
	\end{observation}
	
	\begin{proof}
		This follows immediately from the fact that every edge that contains $v$ also contains $u$.
	\end{proof}
	
	\subsection{Small hypergraphs saturated with respect to pairs}\label{subsection upper bound small graphs}
	
	Our first construction is for a small hypergraph that is Berge-$K_4$-saturated with some nice extra properties.
	
	\begin{construction}\label{construction C for L = 4}
		Fix $k\geq 3$. Let $D=\{d_1,d_2,\dots,d_{k-3}\}$ (when $k=3$ we allow $D=\emptyset$) and $C=\{c_1,c_2,c_3,c_4,c_5\}$ be disjoint sets of vertices. Let $\mC(k,4)$ be the $k$-uniform hypergraph with $V(\mC(k,4))=C\cup D$, and 
		\[
		E(\mC(k,4))=\{\{c_i,c_{i+1},c_{i+2}\}\cup D\mid i\in[5]\},
		\]
		where the indices $i$ are taken modulo $5$.
	\end{construction}
	
	It is worth noting that $\mC(3,4)$ is the $3$-uniform \emph{tight cycle} on $5$ vertices. We now show that $\mC(k,4)$ has the property that every pair is $4$-good and every three vertices form a Berge clique. It is also easy to see that $\mC(k,4)$ is Berge-$K_4$-saturated, but we do not explicitly prove this until Section~\ref{subsection upper bound final construction}.
	
	\begin{lemma}\label{lemma C when L is 4}
		Let $k\geq 3$ and let $\mC:=\mC(k,4)$. Then
		\begin{enumerate}
			\item Every pair in $V(\mC)$ is $4$-good, and\label{lemma C when L is 4 condition good pairs}
			\item Every set of three vertices in $V(\mC)$ are the core vertices of a Berge-$K_3$.\label{lemma C when L is 4 condition berge triangles}
		\end{enumerate}
	\end{lemma}
	
	\begin{proof}
		First, we prove \eqref{lemma C when L is 4 condition good pairs}. Let $xy\subseteq V(\mC)$. First, let us consider the case where $xy=c_ic_j$ for some $i,j\in [5]$. We may assume without loss of generality that $j\in \{i+1,i+2\}$ (modulo $5$). In either case, $\mC+c_ic_j$ contains a Berge-$K_4$ with core vertices $c_{i-1}, c_i, c_{i+1}$ and $c_{i+2}$ with Berge edge map
		\begin{align*}
			c_{i-1}c_i&\mapsto c_{i-2}c_{i-1}c_i\cup D,\\
			c_{i-1}c_{i+1}&\mapsto c_{i-1}c_ic_{i+1}\cup D,\\
			c_{i-1}c_{i+2}&\mapsto c_{i+2}c_{i+3}c_{i-1}\cup D\\
			c_{i+1}c_{i+2}&\mapsto c_{i+1}c_{i+2}c_{i-2}\cup D,\\
			c_ic_{i+1}&\mapsto \begin{cases}
				c_ic_{i+1}&\text{ if }j=i+1,\\
				c_ic_{i+1}c_{i+2}\cup D&\text{ if }j=i+2,
			\end{cases}\\
			c_ic_{i+2}&\mapsto \begin{cases}
				c_ic_{i+1}c_{i+2}\cup D&\text{ if }j=i+1,\\
				c_ic_{i+2}&\text{ if }j=i+2.\\
			\end{cases}
		\end{align*}
		
		Now, if exactly one of $x$ or $y$ is in $D$, say $x\in D$, $y=c_i$, then note that $c_{i+1}\preceq x$, so by Observation~\ref{observation vertex encompass adding an edge}, since $\mC+c_ic_{i+1}$ contains a Berge-$K_4$ with core vertices $c_{i-1}, c_{i}, c_{i+1}$ and $c_{i+2}$, $H+xy$ contains a Berge-$K_4$ as well.
		
		Finally, if both $x,y\in D$, then we note that $c_1\preceq x$ and $c_2\preceq y$, so applying Observation~\ref{observation vertex encompass adding an edge} twice, we first note that since $\mC+c_1c_2$ contains a Berge-$K_4$ with core vertices $c_5$, $c_1$, $c_2$ and $c_3$, $\mC+c_1y$ contains a Berge-$K_4$ with core vertices $c_5$, $c_1$, $y$ and $c_3$, and then $\mC+xy$ contains a Berge-$K_4$.
		
		Now let us focus on \eqref{lemma C when L is 4 condition berge triangles}. Let $xyz\subseteq V(\mC)$. First, we will consider the case when $xyz\subseteq C$. Due to symmetry, we may assume that $xyz=c_1c_2c_3$ or $xyz=c_1c_2c_4$. If $xyz=c_1c_2c_3$, then
		\begin{align*}
			c_1c_2&\mapsto c_5c_1c_2,\\
			c_1c_3&\mapsto c_1c_2c_3,\\
			c_2c_3&\mapsto c_2c_3c_4,
		\end{align*}
		is a Berge edge map for a Berge-$K_3$ with core vertices $x,y,z$. If $xyz=c_1c_2c_4$, then
		\begin{align*}
			c_1c_2&\mapsto c_5c_1c_2,\\
			c_1c_4&\mapsto c_4c_5c_1,\\
			c_2c_4&\mapsto c_2c_3c_4,
		\end{align*}
		again gives us a Berge-$K_3$.
		
		If $xyz\not\subseteq C$, since $c_i\preceq d_j$ for all $i\in [5]$, $j\in [k-3]$, by repeated application of Observation~\ref{observation encompass F to F}, starting with any triple contained in $C$ (that also contains any elements of $xyz$ that are in $C$), we can replace vertices in $C$ with the vertices of $xyz$ that are in $D$, each step retaining the property that the triple is the core of a Berge-$K_3$, resulting in a Berge-$K_3$ with core vertices $xyz$.
	\end{proof}
	
	The following construction and lemma are analogous to Construction~\ref{construction C for L = 4} and Lemma~\ref{lemma C when L is 4}, but for the case $\ell\geq 5$.
	
	\begin{construction}\label{construction C for L = 5+}
		Fix $k\geq 3$, $\ell\geq 5$. Let $D=\{d_1,d_2,\dots,d_{k-3}\}$ and $C=\{c_1,c_2,c_3,\dots,c_{\ell}\}$. Start with a copy of the ($2$-uniform) graph $K:=K_{\ell}-c_1c_2$ on the vertex set $C$, and extend the edges of $K$ to $k$-uniform hyperedges in the following way. Let $e\in E(K)$.
		\begin{enumerate}
			\item If $c_1\in e, c_2\notin e$, then extend $e$ to $e\cup \{c_2\}\cup D$.\label{construction k>4 condition c1 in e}
			\item If $e=\{c_2,c_3\}$, then extend $e$ to $e\cup \{c_4\}\cup D$. \label{construction k>4 condition e is c2c3}
			\item If $e=\{c_2,c_4\}$, then extend $e$ to $e\cup\{c_5\}\cup D$.\label{construction k>4 condition e is c2c4}
			\item If $c_1\notin e,c_2\in e$, and $c_3,c_4\notin e$, then extend $e$ to $e\cup\{c_3\}\cup D$.\label{construction k>4 condition c2 in e c3 c4 not}
			\item If $c_1,c_2\notin e$, then extend $e$ to $e\cup \{c_1\}\cup D$.\label{construction k>4 condition c1c2 not in e}
		\end{enumerate}
		Let $\mC(k,\ell)$ be the resulting $k$-uniform hypergraph.
	\end{construction}
	
	\begin{lemma}\label{lemma C when L is 5+}
		Let $k\geq 3$, $\ell\geq 5$ and let $\mC:=\mC(k,\ell)$. Then
		\begin{enumerate}
			\item Every pair in $V(\mC)$ is $\ell$-good, and\label{lemma C when L is 5+ condition good pairs}
			\item Every set of $\ell-1$ vertices in $V(\mC)$ are the core vertices of a Berge-$K_{\ell-1}$.\label{lemma C when L is 5+ condition berge triangles}
		\end{enumerate}
	\end{lemma}
	
	\begin{proof}
		Let $K$, $C$ and $D$ be defined as in Construction~\ref{construction C for L = 5+}. Notice that the extension of the edges in $E(K)$ to edges in $E(\mC)$ gives a bijection $\phi: E(\mC)\to E(K)$. In order to show that every pair $xy$ in $\mC$ is $\ell$-good, we show that there is a modification of $\phi$ that witnesses a Berge-$K_\ell$ using the edge $xy$. 
		
		\textbf{Case 1:} $xy\subseteq C$. Let $K^*$ be the copy of $K_\ell$ with vertex set $C$. In this case, we will show that $K^*$ is a witness to a Berge-$K_\ell$ in $\mC+xy$. Notice that if $xy=c_1c_2$, then $\phi$ in addition to $xy\mapsto c_1c_2$ gives a Berge-$K_\ell$. 
		
		\textbf{Case 1.1:} $c_1\in xy$, $c_2\not\in xy$, say $x=c_1$. Then by Construction~\ref{construction C for L = 5+}~\eqref{construction k>4 condition c1 in e}, $\phi^{-1}(xy)=c_1yc_2\cup D$, so we let $\phi_{xy}: E(\mC)\cup \{xy\}\to E(K^*)$ be given by 
		\[
		\phi_{xy}(e)=
		\begin{cases}
			xy=c_1y&\text{if }e=xy,\\
			c_1c_2&\text{if  }e=c_1c_2y\cup D,\\
			\phi(e) &\text{otherwise}.
		\end{cases}.
		\]
		
		\textbf{Case 1.2:} $c_2c_3=xy$, say $c_2=x$ and $c_3=y$. We have that $\phi^{-1}(xy)=c_2c_3c_4\cup D$, while $\phi^{-1}(c_3c_4)=c_1c_3c_4\cup D$, and finally $\phi^{-1}(c_1c_3)=c_1c_2c_3\cup D$. Thus, we let $\phi_{xy}: E(\mC)\cup \{xy\}\to E(K^*)$ be given by 
		\[
		\phi_{xy}(e)=
		\begin{cases}
			xy=c_2c_3&\text{if }e=xy,\\
			c_3c_4&\text{if }e=c_2c_3c_4\cup D,\\
			c_1c_3&\text{if }e=c_1c_3c_4\cup D,\\
			c_1c_2&\text{if }e=c_1c_2c_3\cup D,\\
			\phi(e)& \text{otherwise}.
		\end{cases}
		\]
		
		\textbf{Case 1.3:} $c_2c_4=xy$, say $c_2=x$ and $c_4=y$. We have that $\phi^{-1}(xy)=c_2c_4c_5\cup D$, while $\phi^{-1}(c_4c_5)=c_1c_4c_5\cup D$, and finally $\phi^{-1}(c_1c_4)=c_1c_2c_4\cup D$. Thus, we let $\phi_{xy}: E(\mC)\cup \{xy\}\to E(K^*)$ be given by 
		\[
		\phi_{xy}(e)=
		\begin{cases}
			xy=c_2c_4&\text{if }e=xy,\\
			c_4c_5&\text{if }e=c_2c_4c_5\cup D,\\
			c_1c_4&\text{if }e=c_1c_4c_5\cup D,\\
			c_1c_2&\text{if }e=c_1c_2c_4\cup D,\\
			\phi(e)& \text{otherwise}.
		\end{cases}
		\]
		
		\textbf{Case 1.4:} $c_2\in xy$, $c_1,c_3,c_4\not\in xy$, say $c_2=x$. Note that $\phi^{-1}(xy)=c_2yc_3\cup D$, while $\phi^{-1}(c_3y)=c_3yc_1\cup D$, and $\phi^{-1}(c_1y)=c_1yc_2\cup D$. This leads us to $\phi_{xy}: E(\mC)\cup \{xy\}\to E(K^*)$ given by 
		\[
		\phi_{xy}(e)=
		\begin{cases}
			xy=c_2y&\text{if }e=xy\\
			c_3y&\text{if }e=c_2yc_3\cup D\\
			c_1y&\text{if }e=c_3yc_1\cup D\\
			c_1c_2&\text{if }e=c_1yc_2\cup D\\
			\phi(e)&\text{otherwise}
		\end{cases}.
		\]

		\textbf{Case 1.5:}  $c_1c_2\cap xy=\emptyset$. Then $\phi^{-1}(xy)=c_1xy\cup D$, and $\phi^{-1}(c_1x)=c_1xc_2\cup D$. Let $\phi_{xy}: E(\mC)\cup \{xy\}\to E(K^*)$ be given by 
		\[
		\phi_{xy}(e)=
		\begin{cases}
			xy&\text{if }e=xy,\\
			c_1x&\text{if }e=c_1xy\cup D,\\
			c_1c_2&\text{if }e=c_1xc_2\cup D,\\
			\phi(e)&\text{otherwise}.
		\end{cases}
		\]
		
		In all subcases, $\phi_{xy}$ gives a Berge-$K_\ell$.
		
		\textbf{Case 2:}$|xy\cap C|=1$. Assume without loss of generality that $x\in C$ and $y\in D$. let $c\in C\setminus\{x\}$ be any vertex, and note $c\preceq y$. By Case 1, $\mC+xc$ contains a Berge-$K_\ell$ with all core vertices in $C$ (in particular, $y$ is not core). Then by Observation~\ref{observation vertex encompass adding an edge}, $\mC+xy$ contains a Berge-$K_\ell$.
		
		\textbf{Case 3:} $xy\subseteq D$. Note that $c_1\preceq x$ and $c_2\preceq y$. By Case 1, $\mC+c_1c_2$ contains a Berge-$K_\ell$ with all core vertices in $C$. Then by Observation~\ref{observation vertex encompass adding an edge}, $\mC+xc_2$ contains a Berge-$K_\ell$ with all core vertices in $C\cup \{x\}$, and so by Observation~\ref{observation vertex encompass adding an edge} again, $\mC+xy$ contains a Berge-$K_{\ell}$.
		
		Now let us prove statement~\ref{lemma C when L is 5+ condition berge triangles}. of the Lemma~\ref{lemma C when L is 5+}. First, let $X\subseteq C$ be a set of $\ell-1$ vertices, and let $c\in C\setminus X$ and $x\in X$. By Case 1 of statement~\ref{lemma C when L is 5+ condition good pairs}. of Lemma~\ref{lemma C when L is 5+}, $\mC+cx$ contains a Berge-$K_{\ell}$ with all core vertices in $C$. In particular, every vertex in $X$ is core, and $cx\not\subseteq X$, so $\mC$ must contain a Berge-$K_{\ell-1}$ on $X$.
		
		Now assume $X\subseteq V(\mC)$, $|X|=\ell-1$, but with $X\not\subseteq C$. 
		Recall $c_i\preceq d_j$ for all $i\in [\ell]$, $j\in [k-3]$. Therefore, starting with any set $X'$ of $\ell-1$ vertices with $X\cap C\subseteq X'\subseteq C$, we can replace the vertices in $X'\setminus X$ with vertices in $X\cap D$ one at a time by repeated application of Observation~\ref{observation encompass F to F}. 
		At each step, we retain the property that the current set of $\ell-1$ vertices is the core of a Berge-$K_{\ell-1}$, resulting in a Berge-$K_{\ell-1}$ on $X$.
	\end{proof}

	\subsection{Saturated construction for large \texorpdfstring{$n$}{2}}\label{subsection upper bound final construction}
	
	Here we provide the final construction which will provide the desired upper bound in Theorem~\ref{theorem main}.
	
\begin{figure}
		\begin{center}
			\includegraphics[scale=.27]{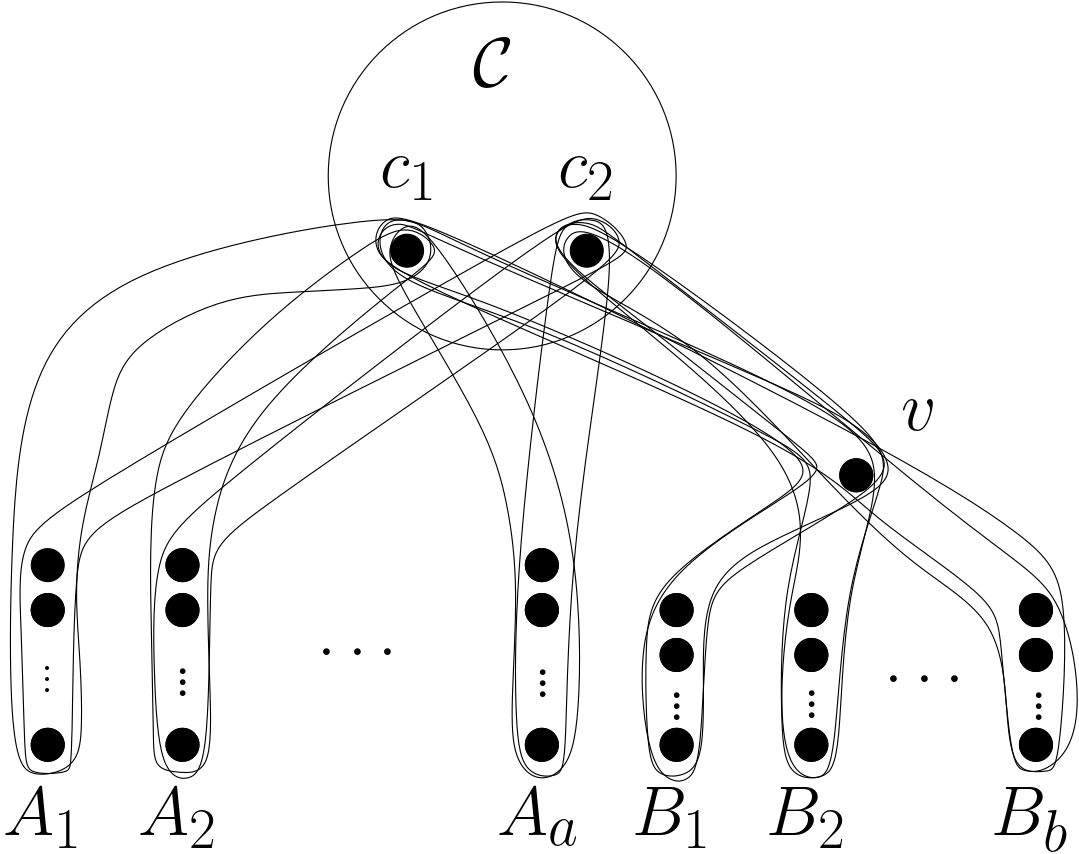}
		\end{center}
		\caption{The construction $\mS(n,k,4)$}\label{figure construction}
	\end{figure}

	\begin{construction}\label{construction final construction}
		Fix $k\geq 3$ and $\ell\geq 4$, and let $n\geq 10k^2\ell$. Let $\mC:=\mC(k,\ell)$, and let $c_1,c_2,\dots,c_{\ell-2}$ be $\ell-2$ distinct vertices from $\mC$. Let $a,b\in \mathbb{Z}_{\geq 0}$ be such that 
		\begin{equation}\label{equation a and b in construction S}
			a(k-1)+b(k-2)=n-|V(\mC)|-1
		\end{equation} 
		and $1\leq b\leq k-1$. We then construct the hypergraph $\mathcal{S}:=\mathcal{S}(n,k,\ell)$ with vertex set 
		\[
		V(\mathcal{S})=V(\mC)\cup \bigcup_{i=1}^a A_i\cup\bigcup_{i=1}^bB_i\cup \{v\}
		\]
		where $|A_i|=k-1$ for all $i\in [a]$ and $|B_i|=k-2$ for all $i\in [b]$, noting that $|V(\mathcal{S})|=|V(\mC)|+a(k-1)+b(k-2)+1=n$. Let 
		\begin{equation}\label{equation edges in construction S}
			E(\mathcal{S})=E(\mC)\cup \{A_i\cup\{c_j\}\mid i\in [a], j\in [\ell-2]\}\cup \{B_i\cup\{v,c_j\}\mid i\in [b], j\in [\ell-2]\}.
		\end{equation}
		See Figure~\ref{figure construction} for a drawing of $\mS(n,k,4)$.
	\end{construction}
	
	The bound $n\geq 10k^2\ell$ is not the best possible, it is simply an easy-to-state bound that is large enough to guarantee that a choice of integers $a,b$ in Construction~\ref{construction final construction} exists.
	
	We now show that $\mS(n,k,\ell)$ is Berge-$K_\ell$-saturated, which will provide the desired upper bound in Theorem~\ref{theorem main}.
	
	\begin{lemma}\label{lemma final construction is saturated}
		The hypergraph $\mathcal{S}:=\mathcal{S}(n,k,\ell)$ is Berge-$K_\ell$-saturated for all $k\geq 3$, $\ell\geq 4$ and $n\geq 10k^2\ell$.
	\end{lemma}
	
	\begin{proof}
		First let us show that $\mS$ is Berge-$K_\ell$-free. Indeed, the only vertices that have degree at least $\ell-1$ are in $V(\mC)\cup \{v\}$. The vertex $v$ only has $\ell-2$ neighbors with degree at least $\ell-1$, so $v$ cannot be a core vertex of any Berge-$K_\ell$, so any Berge-$K_\ell$ would be contained in $V(\mC)$, but only $\binom{\ell}2-1$ edges of $\mS$ have at least two vertices in $V(\mC)$, so no Berge-$K_\ell$ exists in $\mS$.
		
		Now, let $e\in E(\overline{\mS})$. 
		By Lemmas~\ref{lemma C when L is 4}~\eqref{lemma C when L is 4 condition good pairs} and \ref{lemma C when L is 5+}~\eqref{lemma C when L is 5+ condition good pairs}, every pair in $V(\mC)$ is $\ell$-good, so we may assume $|e\cap V(\mC)|\leq 1$.
		
		\textbf{Case 1:} Suppose there exists a pair $X,Y\in \{A_i\mid i\in [a]\}\cup \{B_i\mid i\in [b]\}$ with $X\neq Y$ such that $X\cap e\neq \emptyset$ and $Y\cap e\neq \emptyset$. 
		For simplicity let us assume $X,Y\in \{A_i\mid i\in [a]\}$, as the case where one or more of $X$ or $Y$ are in $\{B_i\mid i\in[b]\}$ is nearly identical. Let $x\in X\cap e$ and $y\in Y\cap e$.
		We claim there is a Berge-$K_\ell$ with core vertices $\{x,y,c_1,c_2,\dots,c_{\ell-2}\}$.
		Indeed, by Lemma~\ref{lemma C when L is 4}~\eqref{lemma C when L is 4 condition berge triangles} and Lemma~\ref{lemma C when L is 5+}~\eqref{lemma C when L is 5+ condition berge triangles}, there is a Berge-$K_{\ell-2}$ with core vertices $\{c_1,c_2,\dots,c_{\ell-2}\}$ using edges from $\mC$ only. 
		Then the edges $X\cup \{c_i\}$ and $Y\cup \{c_i\}$ can play the role of $xc_i$ and $yc_i$ respectively for each $i\in [\ell-2]$ (in the case where $X\in \{B_i\mid i\in [b]\}$, we use the edges $X\cup \{v,c_i\}$), and finally $e$ can play the role of $xy$, completing the Berge-$K_\ell$.
		
		\textbf{Case 2:} Suppose that $e\,\cap (V(\mC)\setminus \{c_1,c_2,\dots,c_{\ell-2}\})\neq\emptyset$. Let $c\in e\,\cap (V(\mC)\setminus \{c_1,c_2,\dots,c_{\ell-2}\})$ and let $x\in e\setminus \{c\}$. Note that $x\not\in V(\mC)$ since $|e\cap V(\mC)|\leq 1$. We will assume $x\in A_1$ since if $x\in A_i$ for $i\in [a]\setminus\{1\}$, $x\in B_i$ for $i\in [b]$ or $x=v$, the case is nearly identical. We claim there is a Berge-$K_{\ell}$ with core vertices in $\{x,c,c_1,c_2,\dots,c_{\ell-2}\}$. Indeed, by Lemma~\ref{lemma C when L is 4}~\eqref{lemma C when L is 4 condition berge triangles} and Lemma~\ref{lemma C when L is 5+}~\eqref{lemma C when L is 5+ condition berge triangles}, there exists a Berge-$K_{\ell-1}$ with core vertices in $\{c,c_1,c_2,\dots,c_{\ell-2}\}$ and only using edges from $\mC$. Then the edge $A_1\cup \{c_i\}$ can play the role of $xc_i$ for $i\in [\ell-2]$, while the edge $e$ plays the role of $xc$, completing the Berge-$K_{\ell}$.
		
		\textbf{Case 3:} Suppose $v\in e$ and neither Case 1 nor Case 2 happen. We claim that there exists a set $X\in \{A_i\mid i\in [a]\}$ such that $e\cap X\neq \emptyset$.
		Indeed, for the sake of contradiction, suppose there does not exist a set $X\in \{A_i\mid i\in [a]\}$ such that $e\cap X\neq \emptyset$.
		Recall that $e$ can contain at most one vertex in $\mC$.
		Since we are not in Case 2, any vertex in $e\cap \mC$ would need to be $c_i$ for some $i\in [\ell-2]$. Since we are not in Case 1, $e$ cannot contain vertices from two $B_j$'s, so all the other vertices in $e$ must come from a single $B_j$, $j\in [b]$, however then $e=\{c_i,v\}\cup B_j\in E(H)$, a contradiction. 
		We may assume, without loss of generality, that $A_1\cap e\neq \emptyset$, say $a\in A_1\cap e$. 
		Then, similar to Case 1, we can find a Berge-$K_{\ell}$ with core vertices in $\{a,v,c_1,c_2,\dots,c_{\ell-2}\}$. 
		Indeed, by Lemma~\ref{lemma C when L is 4}~\eqref{lemma C when L is 4 condition berge triangles} and Lemma~\ref{lemma C when L is 5+}~\eqref{lemma C when L is 5+ condition berge triangles}, there is a Berge-$K_{\ell-2}$ with core vertices $\{c_1,c_2,\dots,c_{\ell-2}\}$ using only edges from $\mC$. 
		Then the edges $A_1\cup \{c_i\}$ and $B_1\cup \{v,c_i\}$ can play the role of $ac_i$ and $vc_i$, respectively for each $i\in [\ell-2]$.
		Finally, $e$ can play the role of $av$, completing the Berge-$K_\ell$.
		
		\textbf{Case 4:} None of the cases 1, 2 or 3 happen.
		We claim that this does not occur. Indeed, $v\not\in e$ since we are not in Case 3. 
		Furthermore, at most one vertex in $\mC$ is in $e$, and if any are, that vertex must be from $\{c_1,c_2,\dots,c_{\ell-2}\}$ since we are not in Case 2.
		Since we are not in Case 1, $e$ intersects at most one set $X\in \{A_i\mid i\in [a]\}\cup \{B_i\mid i\in [b]\}$. 
		Thus, the only way $e$ has $k$ vertices is if $e=X\cup\{c_j\}$ for some $X\in \{A_i\mid i\in [a]\}\cup \{B_i\mid i\in [b]\}$ and $j\in [\ell-2]$. Furthermore, $X$ must be one of the $A_i$'s since the $B_i$'s only have $k-2$ vertices, but $A_i\cup \{c_j\}\in E(\mS)$ for all $i\in [a]$, $j\in [\ell-2]$, contradicting the assumption that $e\in E(\overline{\mS})$.
	\end{proof}
	
	The following theorem summarizes the results of Section~\ref{section upper bound}. Our results could imply a slightly stronger bound than provided below, but we did not attempt to fight for lower-order terms, so this bound suffices for our purposes.
	
	\begin{theorem}\label{theorem upper bound}
		Fix $k\geq 3$, $\ell\geq 4$ and let $n\geq 10k^2\ell$. Then
		\begin{equation}\label{equation upper bound}
			\sat(n,\text{Berge-}K_\ell)\leq \frac{\ell-2}{k-1}n+\binom{\ell}2-1
		\end{equation}
	\end{theorem}
	
	\begin{proof}
		By Lemma~\ref{lemma final construction is saturated}, $\mS:=\mS(n,k,\ell)$ is Berge-$K_\ell$-saturated. Let $\mC:=\mC(k,\ell)$. By \eqref{equation edges in construction S}, we can see that
		\begin{equation}\label{equation first bound on E(S)}
			|E(\mS)|=|E(\mC)|+(a+b)(\ell-2),
		\end{equation}
		where $a$ and $b$ are the integers defined in \eqref{equation a and b in construction S}.
		Recall that $|E(\mC)|=\binom{\ell}2-1$. From \eqref{equation a and b in construction S}, we have that
		\begin{equation}\label{equation bound on a+b}
			a+b=\frac{n-|V(\mC)|-1+b}{k-1}\leq \frac{n}{k-1},
		\end{equation}
		where the inequality follows from the fact that $b\leq k-1$ and the fact that
		\[
		|V(\mC)|=\begin{cases}
			k+2&\text{if }\ell=4,\\
			k+\ell-3&\text{ if }\ell\geq 5.
		\end{cases}
		\]
		Combining \eqref{equation first bound on E(S)} and \eqref{equation bound on a+b} gives us \eqref{equation upper bound}.
	\end{proof}

	\section{Results on linearity for Berge-\texorpdfstring{$F$}{} saturation}\label{section linearity}
	
	In this section, we prove Theorems \ref{theorem linearity delta bound} and \ref{theorem linearity girth bound}. We start with a construction which will help us prove Theorem~\ref{theorem linearity delta bound}.
	
	\begin{construction}\label{4.1}
		Fix $k\geq 3$, let $F$ be a ($2$-uniform) graph with $|V(F)|\geq \alpha(F)+2$, and let $n\geq 10k|V(F)|^3$. Set $\nu:=|V(F)|-\alpha(F)-1$ for convenience. 
		Assume that $k>\nu$, and let $a,t\in \mathbb{Z}_{\geq 0}$ be such that
		\begin{equation}\label{equation linearity defining a}
			a(k-\nu+1)+t=n-\nu
		\end{equation}
		and $0\leq t<k-\nu+1$. We then construct the hypergraph $\mH:=\mH(n,k,F)$ with vertex set
		\[
		V(\mH)=V_1\cup \bigcup_{i=1}^a A_i\cup T,
		\]
		where $|V_1|=\nu$, $|A_i|=k-\nu+1$ for all $i\in [a]$, and $|T|=t$. Let $V_1=\{v_1,v_2,\dots,v_\nu\}$. Let
		\[
		E(\mH)=\{(V_1\cup A_i)\setminus \{v_j\}\mid i\in [a], j\in [\nu]\}.
		\]
	\end{construction}
	
	We do not always expect $\mH(n,k,F)$ to be Berge-$F$-free, but the following lemma shows that if it is Berge-$F$-free, the saturation numbers for Berge-$F$ grow linearly.
	
	\begin{lemma}\label{lemma linearity construction}
		Let $F$ be a graph with $|V(F)|\geq \alpha(F)+2$, $k> |V(F)|-\alpha(F)-1=:\nu$ and $n\geq 10k|V(F)|^3$. If $\mH:=\mH(n,k,F)$ is Berge-$F$-free, then 
		\[
		\sat_k(n,\text{Berge-}F)=O(n).
		\]
	\end{lemma}
	
	\begin{proof}
		Let $\alpha:=\alpha(F)$ and $A:=\bigcup_{i=1}^aA_i$. 
		
		First, we claim that for any set $\{a_1,a_2,\dots,a_{\alpha+1}\}\subset A$ in which $a_i$ and $a_j$ are from different $A_k$'s, there exists a Berge-$(K_\nu\vee \overline{K_{\alpha+1}})$ in $\mH$ with core vertices $\{v_1,v_2,\dots,v_\nu\}\cup \{a_1,a_2,\dots,a_{\alpha+1}\}$ and in which the $v_i$'s correspond to the $K_{\nu}$. Indeed, let us assume without loss of generality that $a_i\in A_i$. Now, for each $i\in [\alpha+1]$ and $j\in [\nu]$, we can use the edge $(V_1\cup A_i)\setminus \{v_{j-1}\}$ to connect $a_i$ to $v_j$ (indices of the $v_j$'s taken modulo $\nu$), giving us a Berge $K_{\nu,\alpha+1}$. Then, from \eqref{equation linearity defining a}, we can see that
		\[
		a=\frac{n-\nu-t}{k-\nu+1}\geq\frac{n-k}{k}>\nu^2+\alpha+1,
		\]
		where in the last inequality, we use our bound on $n$ and that $\nu^2+\alpha+2<10|V(F)|^3$. In particular, this implies that we have at least $\nu^2>\binom{\nu}2$ sets $A_i$ with $i>\alpha+1$. Therefore, there are plenty of edges of the form $(V_1\cup A_i)\setminus \{v_j\}$ for $i\in [a]\setminus [\alpha+1]$ and $j\in [\nu]$ to create the Berge-$K_{\nu}$ on $V_1$. This completes the proof of the claim.
		
		Now, arbitrarily add edges to $\mH$ that do not create a Berge-$F$ until the hypergraph is Berge-$F$-saturated, and call the resulting hypergraph $\mH'$. Let $x\in A$ be a vertex. For the sake of contradiction, assume that
		\[
		d_{\mH'}(x)-d_{\mH}(x)>\binom{|V_1|+t+(\alpha-1)k^2}{k-1}.
		\] 
		This implies that $x$ has $\alpha$ neighbors, which we will denote by $x_1,\dots, x_{\alpha}\in A$ such that 
		\begin{enumerate}
			\item $x_i$ and $x_{i'}$ are not contained in the same $A_j$ for $i\neq i'$, $j\in [a]$, and
			\item there exists an injection $\phi$ from $\{x_1,\dots,x_{\alpha}\}$ into $\{e\in E(\mH')\setminus E(\mH)\mid x\in e\}$ such that $x_i\in \phi(x_i)$ for $i\in [\alpha]$.
		\end{enumerate}
		Indeed, since $d_{\mH'}(x)-d_{\mH}(x)> \binom{|V_1|+t}{k-1}$, $x$ must have a neighbor outside of $V_1\cup T$, say $x_1$, with an edge $e_1\in E(\mH')\setminus E(\mH)$ such that $x,x_1\in e_1$. The edge $e_1$ intersects at most $k$ of the $A_i$'s. Let $B_1$ denote the union of all the $A_i$'s which intersect $e_1$, and note that $|B_1|\leq k|A_1|\leq k^2$. Since $d_{\mH'}(x_1)-d_{\mH}(x_1)> \binom{|V_1|+t+k^2}{k-1}$, $x$ must have a second neighbor, say $x_2$, not in $V_1\cup T\cup B_1$, and let $e_2\in E(\mH')\setminus E(\mH)$ be such that $x,x_2\in e_2$, noting that $e_2\neq e_1$. Then we can define $B_2$ to be the union of all the $A_i$'s that intersect $e_1$ or $e_2$, and note that $|B_2|\leq 2k|A_1|\leq 2k^2$. Continuing inductively, we can find $x_3,\dots,x_{\alpha}$, as claimed.
		
		These edges in $E(\mH')\setminus E(\mH)$, along with the Berge-$(K_\nu\vee \overline{K_{\alpha+1}})$ in $\mH$ with core vertices $\{v_1,v_2,\dots,v_\nu,x_1,x_2,\dots,x_{\alpha},x\}$ give us a Berge-$(K_{\nu+1}\vee \overline{K_{\alpha}})$, which contains a Berge-$F$, a contradiction. Thus, $d_{\mH'}(x)-d_{\mH}(x)\leq \binom{|V_1|+t+(\alpha-1)k^2}{k-1}$ for every vertex $x\in A$.

		We can now wastefully bound $|E(\mH')|$. Indeed, we have that $|E(\mH)|=\nu a\leq \nu n=O(n)$, and 
		\begin{align*}
			|E(\mH')\setminus E(\mH)|&\leq |\{e\in V_1\cup T\}|+|\{e\in E(\mH')\setminus E(\mH)\mid e\cap A\neq\emptyset\}|\\
			&\leq \binom{|V_1|+t}{k}+ak\binom{|V_1|+t+(\alpha-1)k^2}{k-1}\\
			&\leq k\binom{|V_1|+t+(\alpha-1)k^2}{k-1}n=O(n).
		\end{align*}
		Thus,
		\[
		\sat_k(n,\text{Berge-}F)\leq |E(\mH')|= |E(\mH)|+|E(\mH')\setminus E(\mH)|=O(n).
		\]
	\end{proof}
	
	We will need a result from~\cite{AE2019} to deal with an edge case in Theorem~\ref{theorem linearity delta bound}. The result below is significantly weaker than the result in~\cite{AE2019}, but it suffices for our purposes.
	
	\begin{theorem}[\cite{AE2019}]\label{theorem star bound from other paper}
		For any $k\geq 3$, $\ell\geq 1$, 
		\[
		\sat_k(n,\text{Berge-}K_{1,\ell})=O(n).
		\]
	\end{theorem}
	
	We also need a result that handles the case when $k\leq \nu$. Let $\beta(F)$ denote the vertex cover number of $F$.
	
	\begin{theorem}[\cite{EGMT2018}]\label{thm.ksmall} 
		Fix $k\geq 3$. If $F$ is a graph with $\beta(F)\geq k+1$, then $\sat_k(n,\text{Berge-}F)= O(n)$. 
	\end{theorem}
	
	We can now prove the first of our two results on linearity.

	\begin{proof}[Proof of Theorem~\ref{theorem linearity delta bound}]
		First note that if $\alpha(F)=|V(F)|$, then $F$ contains no edges and the saturation function is not defined. Furthermore, if $\alpha(F)=|V(F)|-1$, then $F$ is a star, and by Theorem~\ref{theorem star bound from other paper}, the result follows. Thus, we may assume $|V(F)|\geq \alpha(F)+2$. Suppose that $|V(F)|-\alpha(F)-1=:\nu\geq k$. Then using the fact that $|V(F)|= \alpha(F)+\beta(F)$, we get $\beta(F)-1\geq k$. In this case, Theorem \ref{thm.ksmall} shows that $\sat_k(n,\text{Berge-}F)= O(n)$. Thus, we will assume that $k>\nu$.

		Let $\mH:=\mH(n,k,F)$ as in Construction~\ref{4.1}. We claim that $\mH$ is Berge-$F$-free. Assume to the contrary that there is a copy of $F$ which witnesses Berge-$F$ in $\mH$. We have that $|V(F)\cap A|\geq |V(F)|-\nu=\alpha(F)+1$, and further there must exist an $i\in [a]$ such that $|V(F)\cap A_i|\geq 2$ since otherwise the $\alpha(F)+1$ vertices in $V(F)\cap A$ would form an independent set in $F$ which is too large. Say $x,y\in A_i$. We must have $d_F(x)+d_F(y)-1\leq \nu$ since there are only $\nu$ edges of $\mH$ that intersect $A_i$ and only one can play the role of $xy$, contributing to the degree of both vertices, but $d_F(x)+d_F(y)-1\geq 2\delta(F)-1>|V(F)|-\alpha(F)-1=\nu$, a contradiction.
		
		Thus, $\mH$ is Berge-$F$-free, so by Lemma~\ref{lemma linearity construction}, the result holds. 
	\end{proof}

	We state here a construction from~\cite{EGMT2018} for a $k$-uniform hypergraph on a vertex set $V$ of size $n$ (where $n$ is sufficiently large), which will be used to prove our final result. We will let $f(G)$ denote the vertex feedback number of the graph $G$.
	
	\begin{construction}[\cite{EGMT2018}]\label{construction HnaFS}
		Let $G$ be any graph and let $n$ be sufficiently large. If $f(G)=0$, then for any $a$, let $H_k(n,a,G,\emptyset)$ denote the empty graph on vertex set $V$. If $f(G)\geq 1$, let $S$ be a minimum vertex feedback set of $G$, and let $|E(G[S])|=\ell$. Let the (initially empty) vertex set $V$ be partitioned into three sets $V=V_1\cup V_2\cup V_3$, where $|V_1|=f(G)$, and $|V_2|=(k-2)\ell$. We first add $\ell$ edges between $V_1$ and $V_2$ to create a Berge-$G[S]$ with core vertices in $V_1$. Indeed, arbitrarily label the vertices in $V_1$ with labels from the vertex set of $G[S]$, and then for each $2$-edge $uv$ of $G[S]$, add a $k$-edge that contains the vertices labeled $u$ and $v$ in $V_1$, and $k-2$ vertices in $V_2$ so that after all $\ell$ edges are added, each vertex in $V_2$ has degree $1$.
		
		Now, choose some integer $1\leq a\leq k-1$ such that $a+f(G)\geq k$. If $a$ does not divide $|V_3|$, arbitrarily choose $(|V_3|\mod a)$ vertices, and remove them to form the set $V_3'\subset V_3$, with $|V_3'|=ra$ for some $r\in\mathbb{Z}$. 
		Partition $V_3'$ into $r$ sets of size $a$, and let $\mathcal{M}$ be the collection of these $a$-sets. For each $a$-set $A$ in $\mathcal{M}$, add the $\binom{|V_1|}{k-a}$ edges that contain $A$ and $k-a$ vertices from $V_1$. Call this construction $\mH_k(n,a,G,S)$. 
	\end{construction}

     Our poof of Theorem~\ref{theorem linearity girth bound} will use $a=k-f(G)+1$ or $a=k-1$. However, the original statement of the construction of $\mH_k(n,a,G,S)$ in \cite{EGMT2018} takes $a$ as a more general parameter and facilitates the statement of their theorem.

	\begin{theorem}[\cite{EGMT2018}]\label{theorem linearity paper linear bound}
		Let $F$ be a graph with vertex feedback set $S$, $|S|=f(F)$. Let $1\leq a\leq k-1$ be such that $a+f(F)> k$. If $\mH(n,a,F,S)$ does not contain a Berge-$F$, then 
		\[
		\sat_k(n,\text{Berge-}F)=O(n).
		\]
	\end{theorem}
	
	We are now ready to prove Theorem~\ref{theorem linearity girth bound}
	
	\begin{proof}[Proof of Theorem~\ref{theorem linearity girth bound}]
		Let $S$ be a minimum vertex feedback set of $F$.
        First, consider the case when $f\leq k$ so that $1\leq a=k-f+1\leq  k$. Let $H:=H_k(n,k-f+1,F,S)$. We claim that $H$ is Berge-$F$-free. Indeed, suppose for the sake of contradiction, that $H$ contains a Berge-$F$, and let $F$ be a witness to this Berge-$F$.
		For any $A\in\mathcal{M}$, only $\binom{|V_1|}{k-a}=\binom{f}{f-1}=f<g$ edges of $H$ are incident with vertices in $A$, so no cycle of $F$ can be contained in a set $A$.
		This implies that $V_1$ is a vertex feedback set of $F$, and since $|V_1|=f$, $V_1$ is a minimum vertex feedback set of $F$.

        Similarly, if $f>k$, then we choose $a=k-1$ and let $H:=H_k(n,k-f+1,F,S)$. We claim that $H$ is Berge-$F$-free. Indeed, suppose for the sake of contradiction, that $H$ contains a Berge-$F$, and let $F$ be a witness to this Berge-$F$. For any $A\in\mathcal{M}$, only $\binom{|V_1|}{k-a}=\binom{f}{1}=f<g$ edges of $H$ are incident with vertices in $A$, so no cycle of $F$ can be contained in a set $A$.
		This implies that $V_1$ is a vertex feedback set of $F$, and since $|V_1|=f$, $V_1$ is a minimum vertex feedback set of $F$.
		
		In either case, we proceed with the following argument. For any $v\in V_1$, $v$ must be in a cycle $C$ in $F$ that does not contain any other vertices in $V_1$, since otherwise $V_1\setminus\{v\}$ would be a smaller vertex feedback set. However, this implies that $C$ is contained in $\{v\}\cup A$ for some $A\in \mathcal{M}$, but there are only $f<g$ hyperedges of $H$ that intersect $A$, so there are not enough hyperedges to create a Berge copy of the cycle $C$, and we arrive at a contradiction. 
  
        Thus $H$ is Berge-$F$-free, and by Theorem~\ref{theorem linearity paper linear bound}, the result follows.
	\end{proof}

	\section{Concluding remarks}\label{section conclusion}
	
	The saturation number for Berge-$K_3$ is known exactly, whereas here we were only able to determine the asymptotics of the saturation number for larger cliques. It would be nice to be able to find the exact value, at least in some small cases. The most tractable case is $\sat_3(n,\text{Berge-}K_4)$. By counting edges in $\mS(n,3,4)$ from Construction~\ref{construction final construction} more carefully than is done in Theorem~\ref{theorem upper bound}, we can get the following:
	
	\begin{remark}
		For all $n\geq 10$,
		\[
		\mathrm{sat}_3(n,\text{Berge-}K_4)\leq \begin{cases}
			n&\text{ if }n\text{ is odd,}\\
			n+1&\text{ if }n\text{ is even.}
		\end{cases}
		\]
	\end{remark}
	
	It seems difficult to push the lower bound to match this, but it would be quite interesting to see work in this direction.
	
	\begin{problem}
		Determine $\sat_3(n,\text{Berge-}K_4)$ exactly.
		In particular, is $\sat_3(n,\text{Berge-}K_4)=n$ when $n$ is odd?
	\end{problem}

	\bibliographystyle{amsplain}
	\bibliography{bib}{}
	
\end{document}